\documentclass[review,12pt]{elsarticle}

\usepackage{graphicx}
\usepackage{amssymb}
\usepackage{multirow}
\usepackage{lscape}
\usepackage{xspace}
\usepackage{nicefrac}
\usepackage{booktabs}
\usepackage[colorlinks]{hyperref}
\hypersetup{
  colorlinks=true,
  linkcolor=blue,
  citecolor=blue,
}
\usepackage{orcidlink}

\newcommand{\oldS}{{\sc old}\xspace}
\newcommand{\newS}{{\sc new}\xspace}
\newcommand{\cplex}{{\sc cplex}\xspace}
\newcommand{\gurobi}{{\sc gurobi}\xspace}
\newcommand{\xpress}{{\sc xpress}\xspace}
\newcommand{\mosek}{{\sc mosek}\xspace}
\newcommand{\copt}{{\sc copt}\xspace}
\newcommand{\gcc}{{\sc gcc}\xspace}
\newcommand{\miplib}{{\sc miplib}\xspace}
\newcommand{\lp}{{\sc lp}\xspace}
\newcommand{\milp}{{\sc milp}\xspace}

\newcommand{\parten}{\textsc{par}\oldstylenums{10}\xspace}
\newcommand{\parX}{\textsc{par}\xspace}
\newcommand{\parhun}{\textsc{par}\oldstylenums{100}\xspace}

\usepackage{lineno}

\journal{EURO Journal on Computational Optimization}

\begin{document}

\begin{frontmatter}
\title{Progress in Mathematical Programming Solvers from 2001 to 2020}

\author[TUB,ZIB]{Thorsten Koch\orcidlink{0000-0002-1967-0077}}
\ead{koch@zib.de}
\author[FICO]{Timo Berthold\orcidlink{0000-0002-6320-8154}}
\author[ZIB]{Jaap Pedersen\orcidlink{0000-0003-4047-0042}}
\author[TUB]{Charlie Vanaret\orcidlink{0000-0002-1131-7631}}

\address[TUB]{Technische Universit{\"a}t Berlin,
    Chair of Software and Algorithms for Discrete Optimization, Stra{\ss}e des 17. Juni 135, 10623 Berlin, Germany}
\address[ZIB]{Zuse Institute Berlin, Takustra{\ss}e 7, 14195 Berlin, Germany}
\address[FICO]{Fair Isaac Germany GmbH, Stubenwald-Allee 19, 64625 Bensheim, Germany}

\begin{abstract}
This study investigates the progress made in \lp and \milp solver performance during the last two decades by comparing the solver software from the beginning of the millennium with the codes available today. 
On average, we found out that for solving \lp/\milp, computer hardware got about 20 times faster, and the algorithms improved by a factor of about nine for \lp and around 50 for \milp, which gives a total speed-up of about 180 and 1,000 times, respectively. 
However, these numbers have a very high variance and they considerably underestimate the progress made on the algorithmic side: many problem instances can nowadays be solved within seconds, which the old codes are not able to solve within any reasonable time.
\end{abstract}

\end{frontmatter}


\section{How much did the state of the art in (Mixed-Integer) Linear Programming solvers progress during the last two decades?}

\begin{figure}[htbp!]
\centering
\includegraphics[width=\columnwidth]{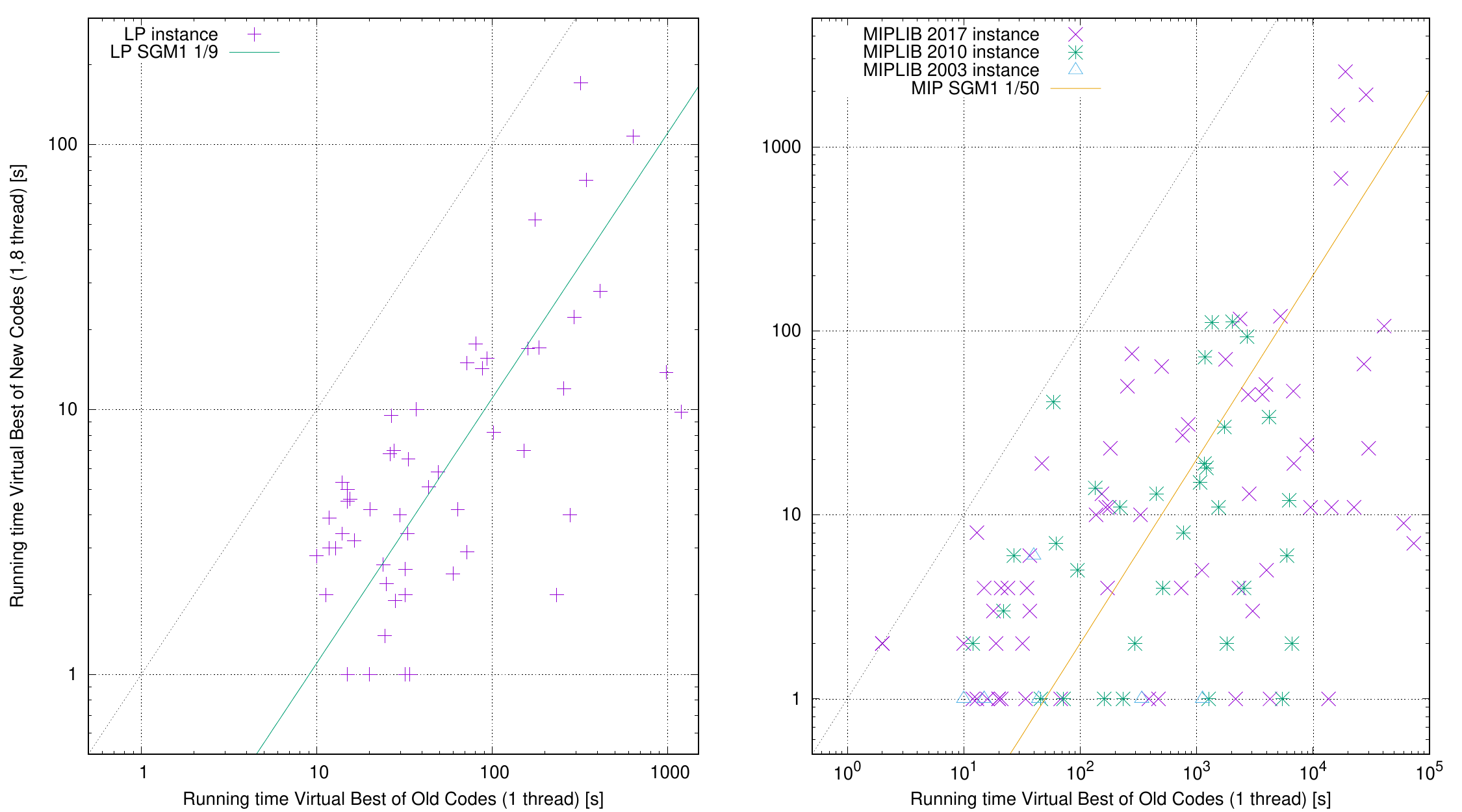}
\caption{Comparison of the running times of various \lp (left) and \milp (right) instances between the virtual best of \cplex~\oldstylenums{7}, \xpress~\oldstylenums{14}, and \mosek~\oldstylenums{3}, from around 2001 and the virtual best of \cplex~\oldstylenums{12.10}, \gurobi~\oldstylenums{9.0}, \xpress~\oldstylenums{8.11}, \mosek~\oldstylenums{8.1}, and \copt~\oldstylenums{1.4} from 2020 running with either 1 or 8 threads on a log scale where \mosek~\oldstylenums{3}, \mosek~\oldstylenums{8.1}, and \copt~\oldstylenums{1.4} are only used on the \lp instances.}
\label{fig:speedup-old-vs-allnew}
\end{figure}

\begin{figure}[htbp!]
\centering
\includegraphics[width=\columnwidth]{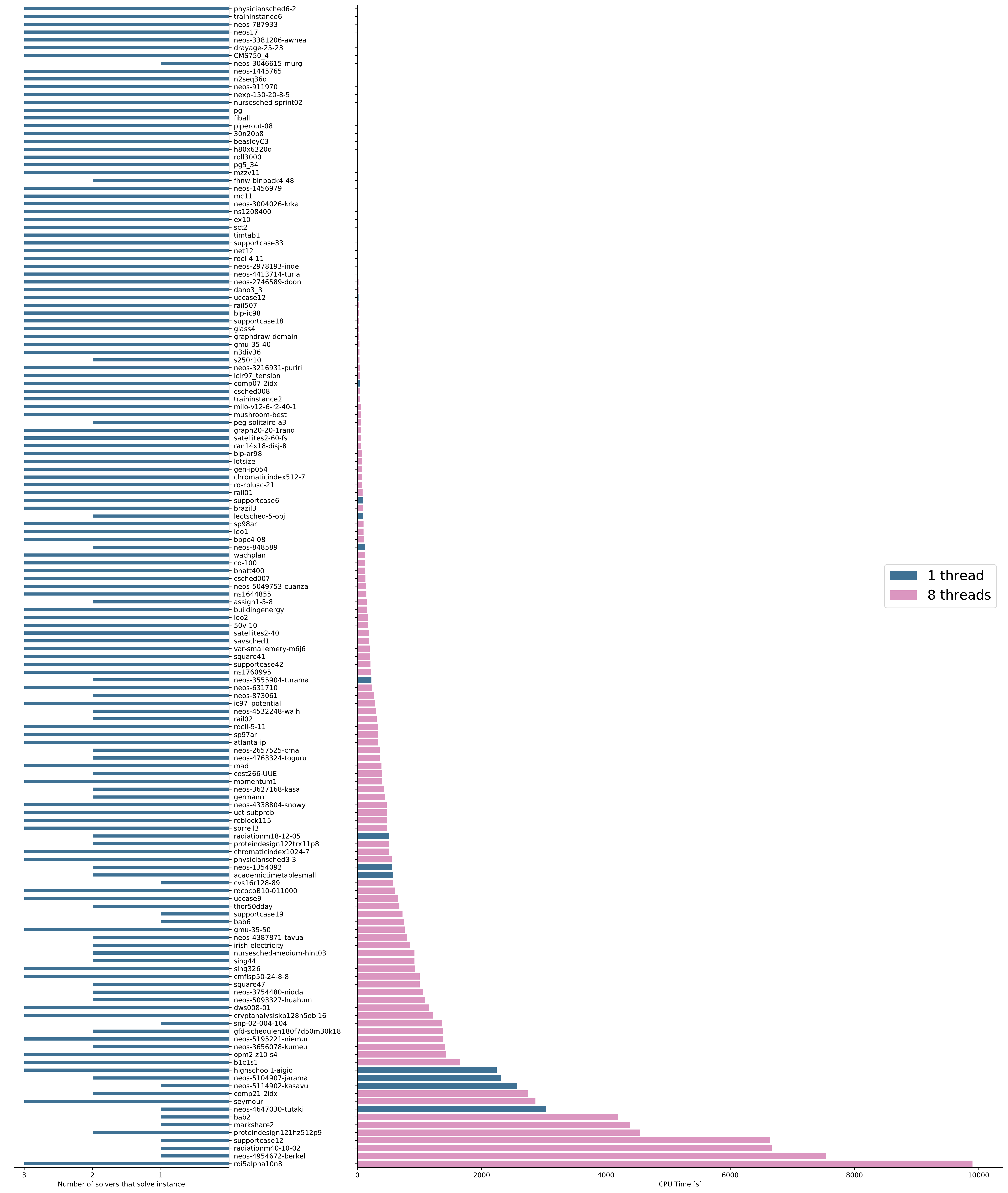}
\caption{Runtime of the virtual best new solver for those 149 instances from the \miplib 2017~\cite{Gleixner2020MIPLIB} benchmark set that could not be solved by any of the old solvers within 24 h. The color for the number of \emph{threads} indicate which was faster. The left column shows how many of the 3 solvers solved the instance within 6 h.}
\label{fig:runtime unsolved part I}
\end{figure}

The present article aims at providing one possible answer to this question. We will argue how progress in \lp and \milp solvers can be measured, how to evaluate this progress computationally, and how to interpret our results. Our findings are summarized in Figures~\ref{fig:speedup-old-vs-allnew} and \ref{fig:runtime unsolved part I}. The main part of this article provides context in which these figures can be interpreted.

Without doubt, computational methods for solving Linear Programs (\lp) and Mixed Integer Linear Programs (\milp) have made tremendous progress during the last 40+ years. The question ``how much?'' naturally arises. And how much of this progress is due to algorithmic improvement compared to advances in hardware and compilers?

\subsection{Previous studies}
This question has been asked before. There are five studies that focus solely on the \cplex solver and cover the 1990s and 2000s. The  first two, by Bixby et al. \cite{BixbyEtal2000,Bixby2002}, investigated the progress from 1987 to 2001 regarding the solution of \lp{}s; the latter concluded: 
\emph{Three orders of magnitude in machine speed and three orders of magnitude in algorithmic speed add up to six orders of magnitude in solving power: A model that might have taken a year to solve 10 years ago can now solve in less than 30 seconds.} For the period from 1997 to 2001, the geometric mean speed-up computed over 677 instances was 2.3. However, it should be noted that the speed-up for large models with more than 500,000 rows was over 20.

Bixby et al. \cite{BixbyEtal2004} examined \milp solving. The study considered 758 instances and compared \cplex 5.0 (released in 1997) and \cplex 8.0 (2002). The geometric mean of the speed-up was about 12.
The speed-up was considerably higher for the instances that required over 27 hours to solve with the older code, reaching an average of 528.

Achterberg and Wunderling \cite{AchterbergWunderling2013} continued the study up to \cplex 12.5 in 2012. The overall geometric mean speed-up on 2,928 \milp models turned out to be 4.71. 
An average speed-up of up to 78.6 was observed for the instances that were challenging for version 8.0.
This is still an underestimation, as the old solver hit the time limit of 10,000 seconds for 732 of the instances, while the new one only had 87 timeouts. 

Lodi \cite{Lodi2010} compared \cplex 1.2 (1991) with \cplex 11.0 (2007) on 1,734 \milp{}s and reported a geometric mean speed-up of 67.9.
Another revealing metric shown is the number of instances solved to optimality within the time limit of 30,000s. On 1,852 \milp{}s, \cplex 1.2 was able to solve a mere 15.0\%, while version 11.0 on the same hardware could solve 67.1\%.

Koch et al.~\cite{KochEtal2011,KochEtal2013} compared the performance of a wide variety of solvers on the \miplib2010. 
The progress from 1996 to 2011 was investigated and the conclusion was unsurprisingly similar. On the one hand, instances that were already solved ``quickly'' did not get solved faster. On the other hand, many instances that used to be  ``difficult'' got solved considerably faster; these were the ones that contributed the most to the overall speed-up. 

Since all of these studies are at least ten years old, it seems about time to give an update on whether \lp and \milp development is still going strong.

\subsection{Setup of this study}
One could argue that all studies, including the present one, have intrinsic biases. The threshold for discarding problems as "too easy" influences the observed speed-up factors. 
The higher the threshold, the higher the speed-up. 
The same happens on the other end: the lower the time limit given to the solver, the lower the achievable speed-up.

Another bias comes from the selection of instances. Instances usually do not enter a collection because they are quickly solved on the first try. Therefore, there is a tendency to collect ``difficult`` instances. On the other hand, modeling practices rely on the efficiency of current solvers, which leads to a selection that under-represents modeling practices that cannot (at the time) be satisfyingly solved. 

Another natural question for our study was which solver to use.
When the initial tests for the \miplib 2010 \cite{KochEtal2011} were performed, all three main commercial solvers achieved roughly the same geometric average running time over the whole benchmark set. The speed difference for individual instances, however, was as large as a factor of 1,000 between the fastest and the slowest solver. Which solver was the fastest was largely instance-dependent. When \miplib 2010 was released, at least one of the three solvers was able to solve each instance within one hour, but it took years until one \emph{single} solver was capable of solving each instance within an hour. To solve a particular instance, why not use the best solver available? Therefore, it seems natural to us that to discuss the overall performance gain, we use the \textit{virtual best} solver available at the time, unless otherwise stated. The term "virtual best" refers to a perfect oracle that would make the right choice among the solvers for a given instance.

In this article, all running times are given for the two virtual solvers \oldS and \newS, where \oldS is the best among 
the \cplex Linear Optimizer 7.0.0 (2000), the Xpress-MP Hyper Integer Barrier Optimizer Release 14.10 (2002), and MOSEK Version 3.2.1.8 (2003) for solving LPs. 
These codes run single-threaded, with the exception of the barrier method \lp solvers within \xpress and \mosek. The best achievable result was systematically kept.
\newS is the best among the \textsc{ibm ilog} \cplex Interactive Optimizer 12.10 (2019), the \gurobi Optimizer 9.0 (2020) and  the \textsc{fico} \xpress Solver 8.11 (2020), as well as \mosek Version 8.1 (2017) and \copt Version 1.4 (2020) for solving \lp{}s.
All solvers were run both sequentially (i.e. single-threaded) and in parallel, allowing the solver to use up to eight threads. 

Our study focuses on the developments of the past twenty years, for three reasons. The first, festive reason is to focus on the period during which EUROPT has been active, following the spirit of this special issue. The second, apparent reason is that this nicely covers the development of \lp and \milp solving in the 21st century (so far). The third, most practical and constraining reason is that it was very tricky to get old and still running solver binaries. As we experienced, a 20-year period is borderline and in some respect already too extensive a time span. There are no contemporary binaries that run on the old 32-bit computers; the old 32-bit binaries failed to run on one of our newer systems. Furthermore, as can be seen in Table~\ref{tab:old-old-vs-new-new}, the speed difference between the old code on the old system and the new code on the new computers is already so enormous that only few instances can be compared in a meaningful way with reasonable effort.

\begin{table}[htbp]
\centering
\sffamily
\footnotesize
\begin{tabular}{l|rr|rr|r}
\toprule
&\multicolumn{2}{c|}{B\&B nodes}&\multicolumn{2}{c|}{Time [s]}\\
Name &  \oldS/P-III & \newS/i7 & \oldS/P-III & \newS/i7 & Speed-up\\
\midrule
nw04    &        131 &      24 &       54 &  1.62 & 33\\
mas76   &    467,454 & 192,719 &      198 &  3.50 & 57\\
neos-1122047&     64 &       1 &       79 &  1.28 & 62\\
mod011  &     16,671 &   3,288 &      926 &  5.77 & 160\\
air05   &      1,961 &   2,189 &      440 &  2.11 & 209\\
qiu     &     36,452 &   4,434 &    1,393 &  1.49 & 935\\
cap6000 &     16,768 &       1 &      268 &  0.10 & 2,680\\
bell5   &    687,056 &     915 &      411 &  0.03 & 13,700\\
neos1171737 & 28,354 &       1 &  116,745 &  2.67 & 43,725\\
\bottomrule
\end{tabular}
\caption{Comparison of selected instances: \oldS (870~MHz Pentium-III) vs \newS (3.6~GHz i7-9700). Remarks: the slowest of \newS with one thread needs 45 s to solve \texttt{mas76} on the i7. In this case the speed-up is exactly the clock ratio between the computers (the solver used in \oldS and the one used in \newS are different though). 
The biggest speed-up happens when the number of B\&B nodes can be reduced. However, whenever there is only one node, no additional speed-up from parallelization occurs.}
\label{tab:old-old-vs-new-new}
\end{table}

\section{Progress in hardware, software and algorithms}
There has been a continuous evolution of the performance of \lp and \milp solvers due to two main, intertwined drivers, namely the development of computers (Section~\ref{fig:progress-computers}) and the algorithmic advances (Section~\ref{fig:progress-algorithms}).
These two sources of progress cannot be easily separated. In the following, we will provide experimental results and discuss which factors influenced the change in performance over time and in which direction.

Unless otherwise stated, all computations have been carried out on an 
8-core, 8-thread Intel Core i7-9700K CPU @ 3.60~GHz with 64~GB of RAM.
It should be noted that modern CPUs adjust their clock speed to the load. The used system might speed up to 
4.7~GHz when running only a single thread. Unfortunately, there is no easy way to track which speed was actually used during a particular run. Therefore, 25\% and more variation of computing time in the measurements are not uncommon. In an experiment, the performance of a single thread halved as we kept the other seven threads busy. For the eight-core runs, the effect is less pronounced, as the machine is already under almost full load by the task to be performed. 

\subsection{Progress in hardware and computational environment}
\label{fig:progress-computers}
In 2001, two of the latest CPUs were the Intel 32-bit Pentium-III at around 1~GHz and the Pentium-4 at 1.5~GHz. IBM offered the 64-bit POWER7 at over 3~GHz. Although 64-bit systems were available, the first 64-bit PC processor was introduced in 2003 and it took quite a few years until more than four gigabytes of memory became standard.
One should keep in mind that there is a gap of several years between the availability of a new architecture and the common use of it by developers and companies. In the following, we list the major developments in hardware and compilers that came into widespread use during the past twenty years:

\begin{description}
\item [Higher clock speed and increased memory bandwidth:] both developments also accelerate old code, even if not recompiled.  
\item [More efficient processing of instructions:] superscalar processors, out-of-order execution, branch prediction, and instruction speed-up. As a consequence, code optimized for an older processor architecture might not perform optimally on a new one. Recompilation is required to fully exploit the improvements.
\item [New instructions:] for example, Fused-Multiply-Add (FMA) and Advanced Vector Extensions (AVX). To exploit these extensions, the code needs at least to be recompiled. The use of highly optimized subroutines (e.g., ATLAS, OpenBLAS, or IMKL) can provide further speed-up. Barrier solvers often have specific subroutines implemented in instruction-set specific assembly code.
\item [Parallel cores and simultaneous multi-threading (SMT):] both have increased the maximal computational performance of a single CPU drastically, however a substantial redesign of the algorithms is required to exploit them. There is almost no automatic benefit for existing codes. Additionally, SMT in particular makes it even harder to determine the best number of parallel threads to use on a given CPU. If memory accesses are the bottleneck, not using SMT can lead to better running times. This is aggravated by the power management of modern CPUs which can decrease clock frequency in case of an increased number of running threads.
\item [Move from 32-bit to 64-bit addressing/processing:] this allows to use more than 4~GB of RAM and to process 64-bit numbers faster. There is no benefit for existing 32-bit codes. Since more memory is used per integer, it possibly can even slow down computations. With some reimplementation however, 64-bit addressing can contribute to performance gains, e.g. by making hash collisions less likely.
\item [Improved optimizing compilers:] from, for example, \gcc version 2.95 (2001) to \gcc version 10 (2020), compilers have improved a lot and generate better performing code. Recompilation is required to benefit from this.
\end{description}

For an overview of hardware and compiler impact on \lp and \milp solver development in the 1980s and 1990s, see~\cite{ashford2007mixed}.

Comparing old and recent architectures is intricate.
Old sequential 32-bit codes using the instruction set available in 2001 will not fully exploit modern architectures. 
Conversely, new parallel 64-bit codes based on recent instructions will not even run on old hardware. 

We performed two small tests to estimate the pure hardware speed-up for solving mathematical optimization problems. First, we ran 33 \lp instances using the single threaded \cplex 7 barrier algorithm without crossover on an old 870~MHZ Pentium-III and the new i7-9700 system, and compared the running times. The speed-up is 21 on average, although it varies between 16 and nearly 47, depending on the particular instance.

Since the requirements of barrier and simplex algorithms are quite diverse, we performed a second test: we solved min-cost-flow problems with the network simplex code in \cplex. 
It can be assumed that this code did not change significantly between version 7.0 and 12.10, as the number of iterations on all instances is identical.
We ran 16 publicly available instances in four different settings:
\cplex 7 on a 870~MHz Pentium-III and on an i7-9700, and we ran \cplex~12 on an otherwise empty i7-9700 and on a fully loaded system.
There is no measurable performance difference between the two \cplex versions regarding the network simplex running on the same hardware. 
\cplex~7 running on the 870~MHZ Pentium-III and on an empty i7-9700 supposedly running at 4.7~GHz boost speed differ by a factor of 20 on average. 
However, if we fully load the system with a multi-core \texttt{stream} benchmark \cite{McAlpinStream}, the performance is halved. One should bear in mind that for each situation, it is not clear where the bottleneck is. The network simplex is known to be highly dependent on the performance of the memory subsystem. 
Overall, the hardware speed-up that we experienced was not constant. The minimum factor is around 15 (i7 empty) and seven (i7 loaded), and the maximum was more than 45.
We would like to point out that small differences in running times are not significant and that the overall impact of the compiler seems small.

The hardware landscape has been changing even more dramatically in the last 10 years, during which Graphics Processing Unit (GPU) accelerators have become widely available. However, as of 2020 (to the best of our knowledge), none of the state-of-the-art \lp/\milp solvers exploits them. Indeed, GPUs are tailored much towards dense processing, while solvers rely heavily on super-sparse linear algebra.

\subsection{Progress in algorithms}
\label{fig:progress-algorithms}

Two decades of research certainly led to significant algorithm improvements. 
One could now ask how much each new feature contributed to the speed-up. Unfortunately, there is no easy and meaningful answer to this. Firstly, we don't know exactly which features were added to each commercial solver. Secondly, since we compare the virtual best, this would be tricky to evaluate even if we knew. Thirdly, as other studies showed \cite{Achterberg2009diss}, for MILP solvers, the sum is more than its parts. In many cases, features support each other. One preprocessing step removing some variables allows another step to remove more. But also the opposite is true: Often, if one component is switched off, part of the effect is provided by the remaining ones. This complicates a meaningful evaluation of feature impact.

The improvements for \milp include many new heuristic methods, such as RINS~\cite{DannaRothbergLePape04} and local branching~\cite{FischettiLodi03}, several  classes of new or improved cutting planes, e.g. MCF cuts~\cite{achterberg2010mcf}, and a large number of \emph{tricks for the bag}, such as conflict analysis~\cite{Achterberg2007}, symmetry detection~\cite{margot2003exploiting}, solution polishing and dynamic search. Most of them either exploit some special structure in the instances, or address a shortcoming in the algorithm for a particular class of problems.
Furthermore, codes have been ported to 64-bit addressing and are therefore able to utilize larger amounts of memory. Moreover, many algorithms have been parallelized~\cite{berthold2018parallelization}, in particular the \milp tree search and barrier methods for \lp solving.

In the area of \lp solving, theoretical progress has been quite limited. There were nonetheless numerous improvements to deal with difficult instances. 
In general, the linear algebra has been sped up by (better) exploiting hyper sparsity and using highly optimized subroutines. 
Preprocessing got better. The parallelization of the barrier methods improved, and there exists nowadays a parallel variant of the simplex algorithm, although its scalability is limited~\cite{HuangfuHall2015}. Nevertheless, with very few exceptions other than due to sheer size, \lp{}s that can be solved nowadays can also be solved with the old codes, provided one is willing to wait long enough. As Figure~\ref{fig:speedup-old-vs-allnew} shows, \lp solvers have become approximately nine times faster since the beginning of the millennium. One should note that often, a given algorithm of a given solver did not become faster, but the fastest choice nowadays is faster than the fastest choice then. Additionally, the ability to solve very large instances has improved considerably.

In the computational study done in 1999, Bixby et al. \cite{BixbyEtal2000} used a 400~MHz P-II with 512~MB of memory. 
The largest \lp that this machine could handle had 1,000,000 rows, 1,685,236 columns, and 3,370,472 non-zeros. 
In a workshop in January 2008 on the \emph{Perspectives in Interior Point Methods for Solving Linear Programs}, the instance \texttt{zib03} with 29,128,799 columns, 19,731,970 rows and 104,422,573 non-zeros was made public. 
As it turned out, the simplex algorithm was not suitable to solve it and barrier methods needed at least about 256~GB of memory, which was not easily available at that time. 
The first to solve it was Christian Bliek in April 2009, running \cplex out-of-core with eight threads 
and converging in 12,035,375 seconds (139 days) to solve the \lp without crossover. Each iteration took 56 hours!
Using modern codes on a machine with 2~TB memory and 4 E7-8880v4 CPUs @ 2.20~GHz with a total of 88 cores, this instance can be solved in 59,432 seconds = 16.5 hours with just 10\% of the available memory used. This is a speed-up of 200 within 10 years. However, when the instance was introduced in 2008, none of the codes was able to solve it. Therefore there was infinite progress in the first year. Furthermore, 2021 was the first time we were able to compute an optimal \emph{basis} solution. 

This pattern is much more pronounced with \milp. First is the step from unsolvable to solved. This is almost always due to algorithmic improvements. Then there is a steady progress both due to algorithmic and hardware improvements until the instance is considered easy. From then on, if any at all, speed-ups are mostly due to hardware only. 
The largest \lp in \xpress' instance collection of practically relevant models has more than 200,000,000 columns and more than 1,300,000,000 non-zeros. It can be solved in about one and a half hours. Solving an instance of this size was impossible with off-the-shelf hardware and software in 2001.
 
Before describing our computational experiments in more detail, note that there are a few caveats to bear in mind:
\begin{itemize}
    \item Since we are interested in the performance of the overall domain, we will compare the virtual best solver \oldS from around 2001  (consisting of \xpress, \cplex, and \mosek) with the virtual best solver \newS from 2020 (consisting of \xpress, \cplex, \gurobi, \mosek, and \copt).

    \item It could be argued that the default parameters are better tuned now and therefore the old codes would benefit more from hand-tuned parameters than the new ones. At the same time, the new codes have more parameters to tune and considerably more sub-algorithms that can be employed. We decided that it is out of the scope of this study to try to hand-tune every instance and therefore only the default values of the solvers will be used.
    
    \item Benchmarking got more prominent and fierce during the last decade, in particular until 2018 (see Mittelmann~\cite{MittelmannBenchmarks}). There has been considerable tuning on the \miplib instances, especially on \miplib 2010 and \miplib 2017. It is fair to say that this clearly benefits the newer solvers and might lead to an overestimation of the progress.
    
    \item
    It should also be noted that instances  
    \texttt{dano3mip, liu, momentum3, protfold} and \texttt{t1717} from \miplib 2003 still remain unsolved, although substantial effort was put into solving them to proven optimality. Furthermore, there are several old instances that still cannot be solved within reasonable time without reformulation or special-purpose codes.
    
\end{itemize}

\section{Computations}

As demonstrated numerous times, the test set and experimental setup have a crucial influence on the outcome of computational studies. The main question is how to handle instances that one actor of the comparison cannot solve. In our case -- not too surprisingly -- \newS is able to solve any instance that \oldS can solve, but not vice versa. When comparing solvers, it is customary to set the maximum run time allowed as time for the comparison. This is reasonable if one compares two solvers on a pre-selected test set.
In our case, the test set and the run time can be chosen; this means that {\bf any speed-up factor can be attained by increasing the number of instances that the old codes cannot solve and increasing the run time allowed. Therefore, we decided to split those questions.}

\subsection{Test set selection}

Our \lp{} test set contains the instances used by Hans Mittelmann for his \lp benchmarks~\cite{MittelmannBenchmarks}, which are listed in Table~\ref{tab:mittelmann-lp-instances}. 

\begin{table}[htbp]
\centering
\footnotesize
\texttt{
\begin{tabular}{llll}
\toprule
buildingenergy  & cont11  & cont1  & cont4  \\ 
datt256  & ds-big  & ex10 & fhnw-binschedule0  \\ 
fome13  & graph40-40  & irish-e & L1\_sixm250obs  \\ 
Linf\_520c  & neos-3025225  & neos-5052403-cy & neos-5251015  \\ 
neos  & nug08-3rd  & pds-100  & physiciansched3-3 \\ 
qap15  & rail02 & rail4284  & rmine15  \\ 
s100 & s250r10 & s82  & savsched1 \\ 
scpm1  & set-cover-model & shs1023 & square41 \\ 
stormG2\_1000  & stp3d  & supportcase10 & tpl-tub-ws1617  \\ 
\bottomrule
\end{tabular}}
\caption{\lp{} instances of Hans Mittelmann's \emph{Benchmark of Simplex LP solvers (1-18-2021)} and \emph{Benchmark of Barrier LP solvers (12-28-2020)}}
\label{tab:mittelmann-lp-instances}
\end{table}
The following instances were excluded from our tests because either \oldS solved them in under 10 seconds or \newS solved them in less than one second: \texttt{chrom1024-7, neos3, ns1688926, self, stat96v1}.

Additionally, we included the following instances from previous Mittelmann \lp benchmarks in the test set: \texttt{dbic1}, \texttt{ns1644855}, \texttt{nug15}.
We further added three larger models, named \texttt{engysys1, engysys2, engysys3}, which we currently use for real-world energy systems research. The motivation for this was to have some hard instances that did not appear in any benchmark so far.

Finally, the \lp{} relaxations of all \miplib-2017 benchmark instances were added to the \lp{} set. Again, we ignored all instances solved by \oldS in under 10 seconds or solved by \newS in less than one second. The resulting instances are listed in Table~\ref{tab:miplib2017-lp-instances}.

\begin{table}[htbp]
\centering
\footnotesize
\texttt{
\begin{tabular}{llll}
\toprule
bab2 & ex9 & highschool1-aigio & k1mushroom \\ 
neos-2075418-temuka & neos-3402454-bohle & neos-5104907-jarama & ns1760995 \\ 
opm2-z10-s4 & rail01 & splice1k1 & square47 \\ 
supportcase19 & triptim1 & & \\
\bottomrule
\end{tabular}}
\caption{\miplib~2017 as \lp{} instances}
\label{tab:miplib2017-lp-instances}
\end{table}

For \lp{}s, \oldS eventually (given enough time) solved all the test instances that \newS could solve, provided the available memory was sufficient. 
One can argue about this selection. Since we compare virtual best solvers, the performance (or lack thereof) of a particular solver does not matter. It is safe to say that these instances were selected because they pose some problems to existing solvers. Therefore, the speed-up experienced on these instances is likely in favor of \newS.

For the \milp comparison we used three versions of \miplib, namely \miplib~2017, \miplib~2010, and \miplib~2003.
It should be noted that we excluded all instances that could not be solved by \newS within six hours, since there is little chance that \oldS would be able to solve them: \texttt{neos-3024952-loue}, \texttt{splice1k1}, \texttt{s100} (from \miplib 2017) and \texttt{ds}, \texttt{timtab2}, \texttt{swath}, \texttt{mkc}, \texttt{t1717}, \texttt{liu}, \texttt{stp3d}, \texttt{dano3mip}, \texttt{momentum3} (from \miplib~2003).
Instances \texttt{mas74} (from \miplib 2017) and \texttt{m100n500k4r1}, \texttt{dfn-gwin-UUM} (from \miplib 2010) were omitted for numerical reasons.
Furthermore, we excluded any instance for which a meaningful speed-up could not be computed due to the limited precision of our timing (whenever \oldS converged within less than ten seconds, and \newS under one second): \texttt{modglob}, \texttt{gesa2}, \texttt{set1ch}, \texttt{irp}, \texttt{p2756}, \texttt{pp08aCUTS}, \texttt{fiber}, \texttt{sp150x300d}, \texttt{neos-827175}, \texttt{drayage-100-23}, \texttt{fixnet6}, \texttt{p200x1188c}, \texttt{gesa2-o}, \texttt{swath1}, \texttt{pp08a}, \texttt{vpm2}, \texttt{neos-1171448}, \texttt{tanglegram2}, \texttt{cbs-cta}. Note that \texttt{nw04} was kept in the test set: although \oldS solved it in just two seconds, \newS needed more than one second (it took two seconds) to solve it. Therefore, a meaningful comparison was possible.

\subsection{Explanation of Figure~\ref{fig:speedup-old-vs-allnew}}
Figure~\ref{fig:speedup-old-vs-allnew} aggregates the results for all instances, both \lp (left) and \milp (right), that could be solved within 24 hours by both \oldS and \newS. In total, these are 56 of 60 \lp and 105 of 339 \milp instances. 
Each symbol represents a single instance. The \emph{x} axis represents the running time of the virtual best old solver and the \emph{y} axis represents the running time of the virtual best new solver. Note that both axes are log-scaled and that we needed one order and two orders of magnitude more to represent the old running times for \lp and \milp instances, respectively. The slowest instance that the old codes could solve in a day took less than three minutes for the new codes. We clipped the times to one second, as this is the precision we could measure.

The dotted grey diagonal is the break-even line. 
Any instance where \oldS is faster than \newS would be above this line. This situation does not occur; however, both virtual solvers took roughly the same time to solve one \milp instance (\texttt{nw04}). Several instances lie on the clipped one-second line, with running times for \oldS of up to 6,000 seconds; all of them have become trivial to solve for \newS.
A reason might be the empirical observation that increasing the allowed running time has a diminishing effect. We will discuss this further in the conclusions related to Figure~\ref{fig:frontier}.

In each plot, a colored diagonal line represents the shifted geometric mean of the speed-up factor for \lp and \milp instances, respectively. All instances (of the corresponding problem type) above (resp. below) the line show a speed-up factor lower (resp. larger) than the mean speed-up. There are some rather extreme cases, in particular for \milp{s}. While the \lp{s} are more concentrated around the mean line, there is no significant difference between \miplib~2010 and \miplib~2017 instances.

The shifted geometric means are one of the main findings. For instances that could be solved by both \oldS and \newS, the pure algorithmic speed-up of \lp solving was about nine in the last twenty years and the speed-up of \milp solving was about 50.

\subsection{Explanation of Figure~\ref{fig:runtime unsolved part I}}
Figure~\ref{fig:runtime unsolved part I} considers the instances that are missing in Figure~\ref{fig:speedup-old-vs-allnew}, that is the 149 (out of 240) instances from \miplib 2017 that could only be solved by new codes. The instances are sorted by the running time of the virtual best new solver. Note that this is a linear scale, not a log scale. Most of the instances (105 of 149) are solved by all the new solvers, whereas twelve instances could only be solved by one of the three solvers within the time limit of 6 h.

Except two instances, all instances were solved in less than one hour and the majority (123 of 149) in less than ten minutes. While one of the eight-threaded solvers was the fastest for most instances, a single-threaded solver won in four out of 13 cases for the hardest instances solved in more than half an hour by \newS.

This seemingly counter-intuitive behavior is explained by the fact that these are instances that are not solved by a tree search, but have a root node that is difficult to solve: the overall solution time is dominated by a hard initial \lp. This explains why they do not benefit from tree-search parallelization and why they still take considerable time with \newS.

Of course, modern \milp solvers also use parallelization at the root node.
However, this mostly happens in terms of concurrent algorithms: different \lp solvers run concurrently to solve the root \lp relaxation; two alternative cut loops run concurrently and the better one carries on at the end of the root node. Concurrent optimization is excellent for hedging against worst-case behavior, but it is inherently slower when the default option would have won in either case. In such a situation, the additional variants (like running primal simplex and barrier for \lp solving or an alternative cut loop) compete for the same resources, and deterministic synchronization might lead to idle times.
In our experiment, the CPU uses turbo-boost for single-thread runs, even amplifying this situation. Thus, we would expect a single-threaded solver to win on an instance that can be solved without much branching and for which dual simplex is the fastest algorithm to solve the initial relaxation.

One main ideal in speeding up the \milp solution process is to reduce the number of branch-and-bound nodes needed. Nearly all modern methods mentioned above (heuristics, cutting planes, conflict analysis, symmetry detection, dynamic search) aim at reducing the number of nodes. The ultimate success is achieved when an instance can be solved in the root node. The progress is visible: among our 339 instances, \oldS solved two instances in the root node, while \newS solved 25 instances.
Unfortunately, the current main direction in hardware development is to increase the number of available threads and the main benefit from parallelization is the ability to process more nodes in parallel. Therefore, hardware and algorithmic development are to a certain extent now non-synergistic. This was different in the past.

\subsection{Take the results with a grain of salt}
We refrained from aggregating instances that could be solved by both \oldS and \newS, and instances that could only be solved by \newS, into a single score. While this might be done by using a time limit for unsolved instances and possibly even a penalty, it can easily skew results. 

\begin{figure}[htbp!]
\centering
\includegraphics[width=0.84\columnwidth]{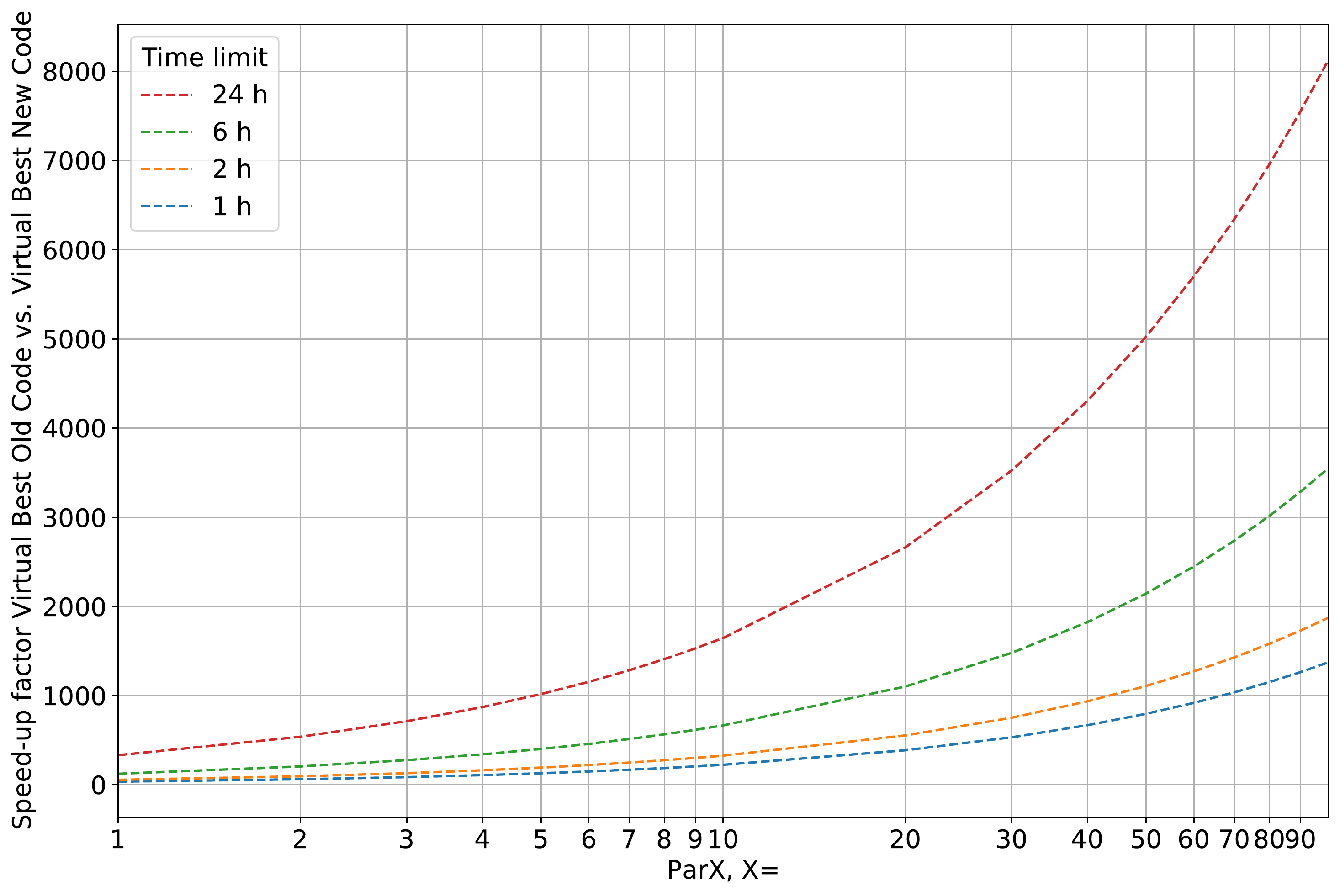}
\caption{Speed-up virtual best \oldS vs. virtual best \newS  on \miplib 2017~\cite{Gleixner2020MIPLIB} benchmark set using different \parX values for various time limits.}
\label{fig:speedup vs parX}
\end{figure}

As a most extreme example, consider using a \parten score (hence weighing all timeouts with a factor of ten times the time limit). Due to the high number of instances that cannot be solved by the old codes, we could obtain almost arbitrarily large speed-up numbers.
With a time limit of one hour and no \parX score (or rather: \parX1), we would get a reasonable speed-up factor of 37 (which is close to the speed-up observed on instances that both versions could solve). Increasing the time limit to 24 hours would give us a speed-up factor of 335. Figure~\ref{fig:speedup vs parX} demonstrates that a similarly large potential for exaggerating results hides in the \parX score. With a time limit of one hour and a \parten score, the ``speed-up factor'' (note that \parten does not actually compute speed-ups) would be 226, with a \parhun score, it would be 1374. Setting the time limit to 24 hours and using \parten, we would get a factor of 1647 and using \parhun, we would get 8125. We see that by driving the time limits up and/or using \parX scores, we can arbitrarily inflate the numbers.

It seems much more sound to report the speed-up factor (50) on solved instances and the impressive 62\% of the instances (156/240) that could be solved by \newS, but not by \oldS. This also shows where the largest progress in \lp and \milp solving lies: 

\centerline{\emph{Making new instances and whole problem classes tractable.}} 

\subsubsection{Performance Variability}

The term performance variability~\cite{lodi2013performance}, loosely speaking, comprises unexpected changes in performance that are triggered by seemingly performance-neutral
changes like changing the order in which constraints of a model are fed into a MILP solver. Besides others, performance variability is caused by imperfect tie-breaking. This results in slight numerical differences caused by floating-point arithmetic which may lead to different decisions being taken during the solution process. Even though one can exploit performance variability in some ways~\cite{Fischetti2015}, it is mainly considered an undesirable property of a solver.

We did the following experiment to investigate whether the amount of performance variability that MILP solvers expose has changed in the past 20 years. Taking only those instances that both \oldS and \newS can solve, we generated ten random permutations of rows and columns for each instance. Mathematically, the instances are equivalent, but different permutations can lead to arbitrarily large differences in run times and tree size, at least for some instances, see, e.g.,~\cite{berthold2018computational}. We ran each permutation of each instance with each solver with a two hour time limit.

We computed the performance variability for \oldS, for \newS using one thread, and for \newS using eight threads. Therefore, we first took the minimum number of nodes needed to solve (or processed within the time limit) a particular instance-permutation combination by either \oldS, by \newS using one thread, or by \newS using eight threads. With these minima, we computed the variability score~\cite{KochEtal2011} of each instance over all permutations, again separately for \oldS, \newS using one thread, and \newS using eight threads. To compensate for seemingly large changes on the lower end, like two nodes instead of one node, which drive up the score tremendously but have no significance in practice, we added a shift of 100 to the node counts.
Finally, we computed the average over the variability scores to get an overall measure of performance variability for \oldS and the two versions of \newS. Table~\ref{tab:variability} summarizes our findings. 

\begin{table}[htbp!]
\centering
\footnotesize
\texttt{
\begin{tabular}{lrrr}
\toprule
    & Mean & Median & Std-dev \\ 
\midrule
\oldS                 & 0.1382 & 0.0926 & 0.1216\\
\newS using 1 thread  & 0.0704 & 0.0577 & 0.0627\\ 
\newS using 8 threads & 0.1137 & 0.0692 & 0.1431\\
\bottomrule
\end{tabular}}
\caption{Mean and median variability scores computed by using ten permutations \label{tab:variability}
}
\end{table}

 We observe that the variability is much lower in the new solvers. The one threaded run of \newS only exposes half of the variability of the one threaded run of \oldS. Unsurprisingly, running in parallel increases variability. 
Note that the standard-deviation is larger than the mean, even though the score is bounded by zero from below. This points to a large spread of variability scores, with a tendency towards the extreme values (including zero variability). The mean is larger than the median because some outliers have considerable variability in all cases.

\section{Conclusion}

\subsection{LP}
For \lp{}s, we computed an average speed-up factor of nine. Combining this with the hardware speed-up, we can conclude that solving \lp{}s got about 180 times faster in the last two decades.
However, the main difference comes from the switch to 64-bit computing, allowing to solve much larger instances, in particular with parallelized barrier codes.
Furthermore, it is fair to say that the solver implementations became ever more refined, leading to extremely stable codes. At the same time, little progress has been made on the theoretical side. 

\subsection{\milp}
In this case, the picture gets more diverse. For the instances solved with the old codes, the average speed-up due to improved algorithms (Figure~\ref{fig:speedup-old-vs-allnew}) is about 50.
This means that solving \milp{}s got faster by 22\% every year during the last 20 years, purely from the algorithmic side, and that is not taking into account how many more instances we can solve today as compared to twenty years ago.
Combining this with the hardware speed-up, we find an average total speed-up for solving \milp{}s of 1,000: fifteen minutes then become less than a second now. 

The most impressive result is shown in Figure~\ref{fig:runtime unsolved part I}: the vastly increased ability to solve \milp instances at all. 149 of 240 instances (62\%) from the \miplib 2017 benchmark cannot be solved by any of the old solvers within a day, even on a modern computer. In contrast, the geometric mean running time for solving these instances by \newS is 104 seconds. We argued why deriving estimated speed-up factors like ``$\nicefrac{1 \textrm{\scriptsize\,day}}{104 \textrm{\scriptsize\,sec}} = 830$-fold speed-up'' would be misleading and one should distinguish between the precise speed-up for instances solved by both and the incredible achievements in solving previously intractable instances.

To summarize, in 2001, one would be pleasantly surprised if one of the solvers would readily solve an arbitrary \milp instance. Nowadays, one is unpleasantly surprised if none of the solvers can tackle it.

\begin{figure}[htbp!]
\centering
\includegraphics[clip=true, trim= 10mm 18mm 120mm 20mm ,width=0.7\columnwidth]{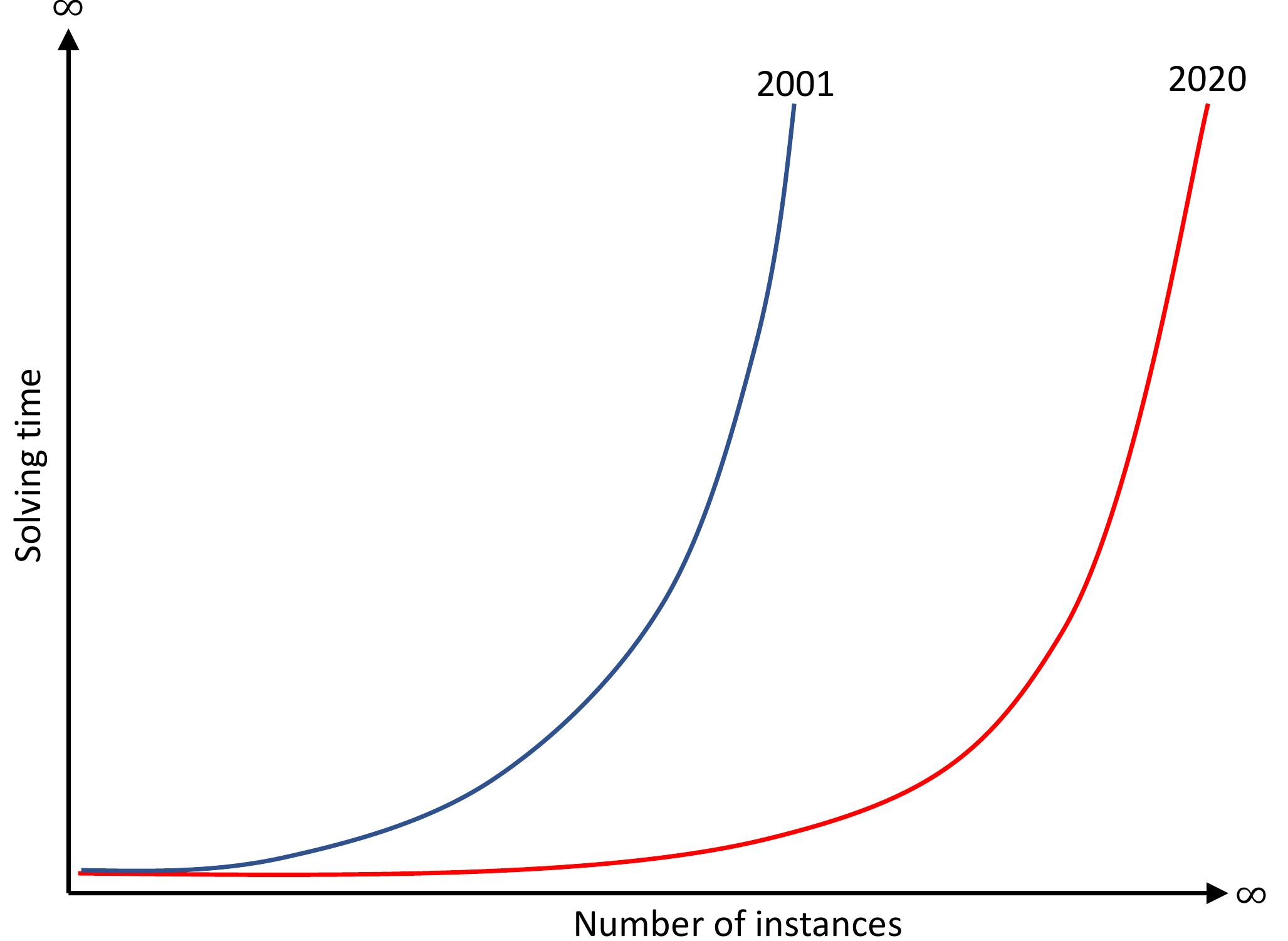}
\caption{The number of instances that are (quickly) solvable is monotonically increasing over time and the frontier of ``difficult'' instances is pushed further to the right (stylized)}
\label{fig:frontier}
\end{figure}

Figure~\ref{fig:frontier} depicts the effect. The number of instances that are solvable right away is ever increasing, but the shape of the frontier stays identical; it is simply pushed to the right. However, it is important to note that the instances on the left are precisely the ones that we wish to solve. 

If we simply speeded up computations, the curve would become more L-shaped but would not be shifted. This is the case with, e.g., (better) parallelization, or if we use the old codes on a new machine. However, it does not change much regarding the overall solvability of instances. To really shift the curve to the right, algorithmic improvements beyond pure computational speed-ups are needed.

\subsection{Outlook}
A problem that we foresee in the future is diminishing returns: as can be deduced from the results, having more and faster cores will not significantly improve the solvability. There are only 14 instances for which \oldS took more than two hours, but less than 24 hours. As \cite{ShinanoAchterbergBertholdetal2020} described, there are individual instances that can be solved by massive amounts of computing power, however there are few of them. A similar situation is true regarding memory. There are, without doubt, some extremely large instances. However, the number of instances that require terabytes of memory are few. And if they do, scaling to higher numbers of cores does not work particularly well because of limited memory bandwidth. There are special algorithms for distributed systems (e.g. \cite{RehfeldtHobbieSchoenheitetal2019}), but these are still far from becoming usable by out-of-the-box solvers. 
Given the change in computer architectures, in particular GPUs, and heterogeneous cores with energy budgets, it becomes increasingly challenging for the solvers to fully exploit the available hardware. This opens interesting directions for research.

A similar observation can be made for the algorithmic side. Overall, it is experienced that every added algorithmic idea affects an increasingly smaller subset of instances.
As always, we hope for breakthrough ideas to appear.
However, so far, solvers still provide significant algorithmic improvements with every release. While additional speed-up by hardware has gone mostly stale, \lp and \milp solvers are still going strong.

\section*{Acknowledgements}
The work for this article has been conducted in the Research Campus MODAL funded by the German Federal Ministry of Education and Research (BMBF) (fund numbers 05M14ZAM, 05M20ZBM).

This work has been supported by the German Federal Ministry of Economic Affairs and Energy (BMWi) through the project UNSEEN (fund no 03EI1004D): \emph{Bewertung der Unsicherheit in linear optimierten Energiesystem-Modellen unter Zuhilfenahme Neuronaler Netze}.

We thank Carsten Dresske for providing us with access to a still nicely running Pentium-III powered computer. 

We thank IBM for providing us with \cplex Version~7, FICO for providing us with \xpress Version~14, and MOSEK for providing us with MOSEK Version~3.

\bibliographystyle{elsarticle-num-names}
\bibliography{biblio.bib}

\end{document}